\input amstex
\documentstyle{amsppt}
\NoBlackBoxes

\mathchardef\t="0074 \mathchardef\o="006F

\mathchardef\p="0070 \mathchardef\r="0072 \mathchardef\t="0074
\mathchardef\b="0062 \mathchardef\d="0064 \mathchardef\i="0069
\mathchardef\n="006E \mathchardef\w="0077 \mathchardef\s="0073
\def\bt{\mathop\boxtimes\limits}
\def\dis{\displaystyle}

\topmatter
\title On spectral gap rigidity \\ and Connes invariant $\chi(M)$
\endtitle
\author SORIN POPA \endauthor

\rightheadtext{On spectral gap rigidity}

\affil University of California, Los Angeles\endaffil

\address Math.Dept., UCLA, LA, CA 90095-155505\endaddress
\email popa\@math.ucla.edu\endemail

\thanks Supported in part by NSF Grant 0601082.\endthanks

\abstract We calculate Connes' invariant $\chi(M)$ for certain
II$_1$ factors $M$ that can be obtained as inductive limits of
subfactors with spectral gap, then use this to answer a question he
posed in 1975, on the structure of McDuff factors $M$ with
$\chi(M)=1$.
\endabstract

\endtopmatter

\document

\heading 1. Introduction \endheading

Given a II$_1$ factor $M$, one denotes by $\chi(M)$ the image in the
outer automorphism group Out$(M)=\text{\rm Aut}(M)/\text{\rm
Int}(M)$, of the group of automorphisms of $M$ that are both
approximately inner and centrally free. This invariant for II$_1$
factors was introduced by Connes in [C75], who used it to solve
several famous problems in von Neumann algebras. On this occasion,
he raised some questions on $\chi(M)$, one of them being whether a
{\it McDuff factor} (i.e. a II$_1$ factor that splits off the
hyperfinite II$_1$ factor $R$) with $\chi(M)=1$ is necessarily of
the form $M=Q\overline{\otimes} R$ with $Q$ a non-Gamma II$_1$
factor, i.e., in the terminology used hereafter, $M$ is {\it
s-McDuff} ({\it strong McDuff})\footnote{I am grateful to V.F.R.
Jones for bringing Connes' question to my attention.}.

We answer this question here, by providing several classes of
examples of McDuff factors who have trivial $\chi(M)$, but are not
s-McDuff. We do this by using the deformation-rigidity methods in
[P01, P03, P06]. Thus, we first show that any infinite tensor
product of non-Gamma II$_1$ factors, $M=\overline{\otimes}_n Q_n$,
satisfies this properties. The second class of examples comes from
subfactor theory. Thus, starting from a non-trivial irreducible
inclusion of non-Gamma factors $P\subset N$ with finite depth, we
consider the enveloping factor $N_\infty$, obtained as the inductive
limit of the associated Jones tower of factors, $P\subset N \subset
N_1 \subset N_2 \subset ....$. We then use [P92] to show that if the
graph $\Gamma=\Gamma_{P,N}$ of $P\subset N$ has the property that
the Perron-Frobenius eigenvector $\vec{s}$ (resp $\vec{t}$) of
$\Gamma^t \Gamma$ (resp. $\Gamma \Gamma^t$) has distinct entries
(e.g. if $\Gamma_{P,N}=A_n$), then $\chi(N_\infty)=1$. On the other
hand, by using deformation-rigidity we prove that $N_\infty$ is not
s-McDuff.

We mention that the deformation of the ambient factors that we use
in our arguments is by {\it inductive limits}, while the rigidity
part is played by the {\it spectral gap} property, considered in
[P06].

\heading 2. Spectral gap and automorphisms in $\text{\rm Ctr}(M),
\overline{\text{\rm Int}(M)}$
\endheading

\noindent {\it 2.1. Definition}. Let $M$ be a II$_1$ factor and
$Q\subset M$ a subfactor. We say that $Q$ has spectral gap in $M$ if
$\forall \varepsilon > 0$, $\exists u_1, ..., u_n\in \Cal U(Q)$ and
$\delta > 0$ such that if $x\in M$ satisfies $\|[x,u_i]\|_2 \leq
\delta \|x\|_2$, $\forall i$, then $\|x-E_{Q'\cap M}(x)\|_2 \leq
\varepsilon \|x\|_2$.

\proclaim{2.2. Lemma} If $Q\subset  M$ is an inclusion of factors
and $Q$ has spectral gap in $M$, then $Q'\cap M^\omega=(Q'\cap
M)^\omega$.
\endproclaim
\vskip .05in \noindent {\it Proof.} Trivial by the definitions.

\hfill $\square$

\vskip .05in

\proclaim{2.3. Lemma} Assume the $\text{\rm II}_1$ factor $M$ is an
inductive limit of subfactors $N_n \nearrow M$ which have spectral
gap in $M$. Then $M'\cap M^\omega=\cap_n (N_n'\cap M)^\omega
=\{(x_n)\mid \exists k_n$ such that $x_n\in N_{k_n}'\cap M$ and
$\lim_\omega k_n=\infty \}$.
\endproclaim
\vskip .05in \noindent {\it Proof.} To see the first equality, note
that we have $M'\cap M^\omega=(\cup_n N_n)'\cap M^\omega=\cap_n
(N_n'\cap M^\omega)$, which by 2.2 is equal to $\cap_n (N_n'\cap
M)^\omega$.

To show the second equality, denote by $\Cal Y$ the given set. We
clearly have $\Cal Y\subset \cap_n (N_n'\cap M)^\omega$. Conversely,
if $(x_n)_n \in \cap_m (N_m'\cap M)^\omega$, then for any $m$ there
exists a neighborhood $V_m$ of $\omega$ such that $\|x_k-E_{N_m'\cap
M}(x_k)\|_2 \leq 2^{-m}$, $\forall k\in V_m$. Moreover, we can take
$V_m\subset V_{m-1}$. Thus, if for each $n\in V_m\setminus V_{m-1}$
we denote $y_n=E_{N_m'\cap M}(x_n)$ and put $k_n=m$, then $(y_n)_n$
coincides with $(x_n)_n$ in $M^\omega$ and clearly $y\in \Cal Y$.
\hfill $\square$

\proclaim{2.4. Lemma} Assume that $Q\subset P$ is an inclusion of
factors. Let $0< \varepsilon < 1/2$ and set $\delta= (\varepsilon/6)^8$.
If $\theta\in \text{\rm Aut}(P)$ satisfies
$\|\theta(v)-v\|_2\leq \delta$, $\forall v\in \Cal U(Q)$, then there
exists a partial isometry $w \in \Cal U(P)$ such that
$\theta(x)w=wx$, $\forall x\in Q$, and
$\|w-1\|_2 \leq \varepsilon$. If in addition $Q'\cap P$ is a
factor, then there exists $u\in \Cal U(P)$ such that $u$ still satisfies
$\theta(x)u=ux$, $\forall x\in Q$, and $\|u-1\|_2 \leq 2\varepsilon$.
\endproclaim
\vskip .05in \noindent {\it Proof.} Let $K=\overline{\text{\rm
co}}^w \{\theta(v)v^*\mid \Cal U(Q)\}$. Note that $\|y\|\leq 1$ and
$\|y-1\|_2\leq \delta$, for all $y\in K$.
Since $K$ is convex and weakly compact in $P\subset L^2(P)$, there
exists a unique element $\xi$ of minimal norm $\|\cdot \|_2$ in $K$.
Since $\theta(v)\xi v^*\in K$ and $\|\theta(v)\xi
v^*\|_2=\|\xi\|_2$, $\forall v\in \Cal U(Q)$, by the uniqueness of
$\xi$ it follows that $\theta(x)\xi = \xi x$, $\forall x\in M$.
Thus, the partial isometry $w\in M$ in the polar decomposition of
$\xi$ still satisfies $\theta(x)w = w x$, $\forall
x\in M$. Moreover, by [C76] we have $\|w-1\|_2 \leq
6\delta^{1/8}=\varepsilon$.

If in addition $Q'\cap P$ is a factor and we assume $\varepsilon
< 1/2$, then $\tau(ww^*)\geq 1/2$ and we can
take partial isometries $e_{12}\in Q'\cap M$, $f_{21}\in
\theta(Q)'\cap M$ such that
$e_{12}e_{12}^*\leq w^*w$, $e_{12}^*e_{12} + ww^*=1$,
$f_{21}^*f_{21}\leq ww^*$, $f_{21}f_{21}^* + ww^*=1$.
Thus, if we define
$u=w + \theta(f_{j1})we_{1j}$, then $u$ is a unitary
satisfying $\theta(x)u=ux$, $\forall x\in Q$,
and $\|u-1\|_2 \leq 2\varepsilon$, proving the statement.
\hfill $\square$

\proclaim{2.5. Lemma} Let $\theta\in \text{\rm Aut}(M)$. Then we
have:

$(a)$ If $N_0\subset M$ has spectral gap and $\theta(x)=\lim_n
u_nxu_n^*$, $\forall x\in M$, then
$$
\lim_{n,m\rightarrow \infty} \|u_{n}^*u_m - E_{N_0'\cap
M}(u_n^*u_m)\|_2=0, \forall k.$$ If in addition $\theta_{|N_0}=id$,
then there exist $v_n\in \Cal U(N_0'\cap M)$ such that $\theta(x)=
\lim_n v_n x v_n^*$, $\forall x\in M$.

$(b)$ Assume $M$ is an inductive limit of subfactors $N_n \nearrow
M$ with spectral gap, such that $N_k'\cap M$ is a factor, $\forall
k\geq 0$. Then $\theta\in \overline{\text{\rm Int}(M)}$ iff $\exists
v_k\in \Cal U(N)$ such that $\theta_{|N_k}=\text{\rm
Ad}(v_k)_{|N_k}$, $\forall k\geq 0$. Moreover, if this is the case
then, $U_n=v_0^*v_n\in N_0'\cap M$ and we have $\theta(x) =
\text{\rm Ad} v_0 \lim_n \text{\rm Ad}U_n (x)$, $\forall x\in M$.
\endproclaim
\vskip .05in \noindent {\it Proof.} The first part of $(a)$ is
trivial. To prove part the second part, note that the spectral gap
of $N_0$ in $M$ and the condition $\lim_n u_nxu_n = x$, $\forall
x\in N_0,$ imply $\lim_n \|u_n-E_{N_0'\cap M}(u_n)\|_2 = 0$. By
[C76], there exist unitary elements $v_n\in \Cal U(N_0'\cap M)$ such
that $\lim_n \|v_n-E_{N_0'\cap M}(u_n)\|_2 = 0$, and thus $\lim_n
\|u_n-v_n\|_2=0$ as well, showing that $\theta=\lim \text{\rm
Ad}(v_n)$.

To prove part $(b)$, assume $\theta=\lim_n \text{\rm Ad}(u_n)$, fix
$k\geq 1$ and let $n_0$ be such that $\|u_{n}^*u_m - E_{N_k'\cap
M}(u_n^*u_m)\|_2\leq 1/4$, $\forall n,m\geq n_0$ (by $(a)$). Then,
for $x\in \Cal U(N_k)$ we have $u_{n_0}^*\theta(x)u_{n_0} =
\lim_{m\rightarrow \infty} u_{n_0}^*u_mxu_m^*u_{n_0}$ and thus

$$
\|u_{n_0}^*\theta(x)u_{n_0} - x\|_2 \leq 2\limsup_m \|u_{n_0}^*u_m -
E_{N_k'\cap M}(u_{n_0}^*u_m)\|_2\leq 1/2 $$

But this implies $\text{\rm Ad}(u_{n_0}) \circ \theta$ is inner on
$N_k$, by Lemma 2.4. Thus, there exists $V_k\in \Cal U(M)$
such that $v_k=u_{n_0}V_k$ satisfies $\theta(x)=\text{\rm Ad}v_k$,
$\forall x\in N_k$.
 \hfill $\square$

\proclaim{2.6. Lemma} Assume that $N_n'\cap M$ is factor, $\forall
n$. Then $\theta\in \text{\rm Ctr}(M)$ iff there exists $n$ and
$u\in \Cal U(M)$ such that $\theta_{|N_n'\cap M} = \text{\rm
Ad}(u)$.
\endproclaim
\vskip .05in \noindent {\it Proof.} It is trivial to see that if
there exists $n$ and $u\in \Cal U(M)$ such that $\theta(x)=uxu^*$,
$\forall x\in N_n'\cap M$, then $\theta\in \text{\rm Ctr}(M)$.

To prove the converse, by Lemma 2.4 it is sufficient to show that
there exists $n$ such that $\|\theta(v)-v\|_2 \leq 1/2$, $\forall
v\in \Cal U(N_n'\cap M)$. Assume on the contrary that $\forall n$,
$\exists v_n \in \Cal U(N_n'\cap M)$ such that
$\|\theta(v_n)-v_n\|_2 > 1/2$. But this implies that $v=(v_n)\in
M'\cap M^\omega$ satisfies $\theta(v)\neq v$, contradicting the
central triviality of $\theta$. \hfill $\square$

\proclaim{2.7. Corollary} Assume the $\text{\rm II}_1$ factor $M$ is
an inductive limit of subfactors $N_n \nearrow M$ such that $N_n$
has spectral gap in $M$ and $N_n'\cap M$ is a factor, $\forall n$.

\vskip .05in \noindent $1^\circ$ If $\theta\in \text{\rm Ctr}(M)\cap
\overline{\text{\rm Int}(M)}$ then for any large enough $n$, there
exists $u,v\in \Cal U(M)$ such that $\theta=\text{\rm Ad}u$ on $N_n$
and $\theta=\text{\rm Ad}v$ on $N_n'\cap M$. Moreover, if for such
an $n$ we denote $\theta'=\text{\rm Ad}u^* \theta$, for any
$\varepsilon >0$, there exists $m\geq n$ and a non-zero partial
isometry $w\in N_n'\cap M$ such that $\theta(x)w=wx$, $\forall x\in
N_m'\cap M$, and $\|w-1\|_2 \leq \varepsilon$.

\vskip .05in \noindent $2^\circ$ If in addition $\Cal N_M(N_n'\cap
M)$ acts innerly on $N_n'\cap M$, then $\text{\rm Ctr}(M)\cap
\overline{\text{\rm Int}(M)}/\text{\rm Int}(M)$ naturally embeds
into the group of automorphisms $\theta\in \text{\rm Aut}(M)$ which
act trivially on $N_n\vee N_n'\cap M$.
\endproclaim
\vskip .05in \noindent {\it Proof.} By Lemmas 2.5 and 2.6, for any
large enough $n$ there exists $u, v\in \Cal U(M)$ such that
$\theta'=\text{\rm Ad}(u)\theta$ is trivial on $N_n$ and equal to
Ad$(v)$ on $N_n'\cap M$. The rest of $1^\circ$ follows from Lemma
2.4.

Under the additional assumption in $2^\circ$, it follows that
$\theta$ can be perturbed by an element in $\text{\rm Int}(N_n'\cap
M)$ so that to act trivially on $N_n\vee N_n'\cap M$. \hfill
$\square$

\heading 3. Deformation-rigidity lemma and examples
\endheading

\proclaim{3.1. Lemma} Let $Q\subset M$ be an inclusion of factors
and assume that $Q$ has spectral gap in $M$. If $N_n\subset M$ are
von Neumann subalgebras such that $\lim_n \|E_{N_n}(x)-x\|_2=0$,
$\forall x\in M$, then for any $\varepsilon > 0$, there exists $n$
such that $N_n'\cap M \subset_{\varepsilon} Q'\cap M$.
\endproclaim
\vskip .05in \noindent {\it Proof.} Since $Q$ has spectral gap in
$M$, by definition there exist $\delta>0$ and $u_1, ..., u_m\in \Cal
U(Q)$ such that if $x\in (M)_1$ satisfies $\|[u_i, x]\|_2 < \delta$,
$\forall i$, then $x\in_\varepsilon Q'\cap M$.

By the hypothesis, there exists $n$ such that
$\|E_{N_n}(u_i)-u_i\|_2 \leq \delta/2$, $\forall i$. Thus, if $x\in
N_n'\cap M$ then

$$
\|[x, u_i]\|_2 \leq 2\|x\| \|E_{N_n}(u_i)-u_i\|_2 \leq \delta,
\forall i,
$$
implying that $x\in_\varepsilon Q'\cap M$. \hfill $\square$

\vskip .05in

We next look for conditions which are sufficient for the assumptions
in Corollary 2.7 to be satisfied and for which Lemma 3.1 thus
applies. As in [P94], we denote by $\Cal G_{P,N}$ the standard
invariant of an inclusion of factors with finite index $P\subset N$
and by $\Gamma_{P,N}$ its standard graph.

\proclaim{3.2. Proposition} $1^\circ$ If $N$ is a non-Gamma
$\text{\rm II}_1$ factor and $S$ is an arbitrary finite factor, then
$N$ has spectral gap in $M=N\overline{\otimes} S$. Moreover, any
subfactor of finite index $P\subset N$ has spectral gap in $M$.

\vskip .05in \noindent $2^\circ$ If $N$ has the property $\text{\rm
(T)}$ and $M$ is a $\text{\rm II}_1$ factor containing $M$, then $N$
has spectral gap in $M$. In particular, if $P$ is a subfactor of
finite index of a property $\text{\rm (T)}$ factor $N$, $P\subset N
\subset N_1 \subset ...$ is the Jones tower and $N_\infty$ the
associated enveloping algebra, then $N_n$ has spectral gap in
$N_\infty$, $\forall n$.

\vskip .05in \noindent $3^\circ$ If $N$ is non-Gamma, $P\subset N$
is a subfactor with finite depth, $N_n \nearrow N_\infty$ the
associated tower and enveloping algebra, then $N_n$ has spectral gap
in $N_\infty$, $\forall n$.

\vskip .05in \noindent $4^\circ$ Let $\Cal G$ be a standard
$\lambda$-lattice, $Q$ a $\text{\rm II}_1$ factor and denote
$P=M_{-1}^{\Cal G}(Q)\subset M_0^{\Cal G}(Q)=N$ the inclusion of
$\text{\rm II}_1$ factors with standard invariant $\Cal G$, as
constructed in ${\text{\rm [P94], [P98]}}$, with $N_n \nearrow
N_\infty$ the associated tower and enveloping algebra. Then $N_n$
has spectral gap in $N_\infty$, $\forall n$. Also, if $Q=L(\Bbb
F_\infty)$ then $P\simeq N_n \simeq L(\Bbb F_\infty)$, $\forall n$.

\vskip .05in

Moreover, in case $3^\circ$, $N_n'\cap N_\infty$ is a factor,
$\forall n$. In turn, in examples $2^\circ$ and $4^\circ$, $N'\cap
N_\infty$ is a factor if and only if the standard graph of $P\subset
N$ is ergodic. This is the case if for instance $\Cal G_{N,M}$ is
strongly amenable, or if $P\subset N$ has graph $A_\infty$, i.e.
when $[N:P]\geq 4$ and the relative commutants $N'\cap N_n$ are
generated by the Jones projections $($i.e. $\Cal G_{P,N}$ is the
so-called Temperley-Lieb standard lattice$)$.

\endproclaim
\vskip .05in \noindent {\it Proof.} Part $1^\circ$ is essentially
due to Connes (see [C76]) and $2^\circ$ is trivial. Part $4^\circ$
is immediate by the proofs of (7.1 and 7.3 in [P90]).

To prove part $3^\circ$, note first that  if $N\supset P \supset P_1
\supset P_2...$ is a tunnel, then $P_k$ has spectral gap in $N$, and
thus in $N\vee N'\cap N_\infty$, for any $k$. Since $P\subset N$ has
finite depth, so does $N\vee N'\cap N_\infty \subset N_\infty$ (see
e.g. [Oc87] or [EK98]) and there exists $k$ such that $P_k'\cap
N_\infty$ contains an orthonormal basis $\{m_j\}_j$ of $N_\infty$
over $N\vee N'\cap N_\infty$. Thus, if $\xi\in N'\cap
L^2(N_\infty)^\omega$, then $\xi=\Sigma_j m_j \xi_j$, for some
unique ``coefficients'' $\xi_j$ lying in $L^2(N\vee N'\cap
N_\infty)^\omega$. Since $\xi$ and $m_j$ commute with $P_k$, it
follows that $\xi_j$ commute with $P_k$ as well, so in fact all
$\xi_j$ lie in $(P_k'\cap N)\vee L^2(N'\cap N_\infty)^\omega$.
Altogether, $\xi\in N'\cap L^2(P_k'\cap N_\infty)^\omega =
L^2(N'\cap N_\infty)^\omega$.

\hfill $\square$

\heading 4. Calculations of $\chi(M)$ and an answer to Connes'
question
\endheading

\proclaim{4.1. Theorem}  If $M=\overline{\otimes}_k Q_k$, with $Q_k$
a sequence of non-Gamma $\text{\rm II}_1$ factors, then $\text{\rm
Ctr}(M)\cap \overline{\text{\rm Int}(M)}=\text{\rm Int}(M)$
$($equivalently $\chi(M)=1)$ and $M$ is McDuff but not s-McDuff.
\endproclaim
\vskip .05in \noindent {\it Proof}. By Corollary 2.7.2$^\circ$ we
have $\chi(M)=1$ (this calculation was in fact already done in
[C75]). Assume $M=Q\overline{\otimes} R$ for some non-Gamma factor
$Q$. Thus, $Q$ has spectral gap in $M$. Denote
$N_n=\overline{\otimes}_{k\leq n} Q_k$. By Lemma 3.2, there exists
$n$ such that $\overline{\otimes}_{k>n} Q_k=N_n'\cap M
\subset_{\varepsilon} Q'\cap M=R$. By [P03], this implies there
exists an isomorphism of $\overline{\otimes}_{k>n} Q_k$, which is a
non-amenable factor, into an amplification of $R$. Since $R$ is
amenable, this is a contradiction. \hfill $\square$

\vskip .05in

From here on, we consider the following special case of inductive
limits of factors: We let $P\subset N$ be a subfactor of finite
Jones index, $P\subset N \subset N_1 \subset N_2 \subset ... $ its
Jones tower and $M=N_\infty=\overline{\cup_n N_n}$ the associated
{\it enveloping} II$_1$ {\it factor}. Under the assumption that $N$
has spectral gap in $N_\infty$, we relate Connes' $\chi$-invariant
of $N_\infty$ with Kawahigashi's $\chi$-invariant of the inclusion
$N'\cap N_\infty\subset P'\cap N_\infty$ and use this to calculate
$\chi(N_\infty)$ for many enveloping factors. In particular, this
will provide more examples of factors with trivial $\chi$-invariant
which are McDuff but not s-McDuff.

Thus, recall from [K93] that if $S\subset R$ is an inclusion of
finite von Neumann algebras then

$$
\chi(R,S) \overset \text{def} \to = \overline{\text{\rm Int}(R,S)}
\cap \text{\rm Ctr}(R,S)/ \text{\rm Int}(R,S),
$$
where $\text{\rm Int}(R,S)$ is the group of inner automorphisms of
$R$ implemented by unitaries in $S$, $\overline{\text{\rm
Int}(R,S)}$ its closure and $\text{\rm Ctr}(R,S)$ the group of
automorphisms $R$ leaving $S$ invariant and acting trivially on
$R'\cap S^\omega$. We consider this invariant for inclusions of the
form $S=N'\cap N_\infty \subset P'\cap N_\infty=R$, where $P\subset
N$ is a (proper) subfactor of finite index and $N_\infty$ is its
enveloping algebra, as above.

\proclaim{4.2. Theorem}  If $N$ has spectral gap in $N_\infty$ then
there is natural isomorphism $\chi(N_\infty)\simeq \chi(P'\cap
N_\infty, N'\cap N_\infty)$.
\endproclaim
\vskip .05in \noindent {\it Proof.} Note first that by [PP83], if
$\sigma$ is an automorphism of $R=P'\cap N_\infty$ leaving $S=N'\cap
N_\infty$ invariant, then there exists $u\in \Cal U(N'\cap
N_\infty)$ such that $u \sigma (e_1) u^*=e_1$, where $e_1=e_P$.
Thus, any element in $\overline{\text{\rm Int}(R,S)}/\text{\rm
Int}(R,S)$ can be represented by an automorphism of the form $\sigma
= \lim_n \text{\rm Ad}(u_n)$, with $u_n \in \Cal U(N_1'\cap
N_\infty)$. Denote by $\Cal G\subset \text{\rm Aut}(R,S)$ the group
of automorphisms of this form. Denote also $\Cal G_0=\Cal G\cap
\text{\rm Ctr}(R,S)$.

Let also $\Cal H\subset \text{\rm Aut}(N_\infty)$ be the group of
automorphisms in $\overline{\text{\rm Int}(N_\infty)}$ which act
trivially on $N_1$ and $\Cal H_0=\Cal H\cap \text{\rm
Ctr}(N_\infty)$.

If $\sigma\in \Cal G$, then let $\psi=\Psi(\sigma)$ be the Banach
limit, $\psi(x)\overset \text{def} \to =\text{\rm Lim}_n u_nxu_n^*$,
$x\in N_\infty$. Thus $\psi$ is a unital trace preserving c.p. map
on $N_\infty$ and by the choice of $u_n$, we have $\psi(x)=x$ for
$x\in N$, $\psi(e_1)=e_1$ and $\psi(x)=\sigma(x)$, $x\in P'\cap
N_\infty$. Since a unital trace preserving c.p. map is
multiplicative on the space of elements on which it is
$\|\cdot\|_2$-isometric, and since the algebra generated by $N$,
$e_1$ and $N'\cap N_\infty$ is $\|\cdot\|_2$-dense in $N_\infty$, it
follows that $\psi\in \overline{\text{\rm Int}(N_\infty)}$.
Moreover, since $N$ has spectral gap in $N_\infty$, $N_\infty'\cap
N_\infty^\omega=N_\infty'\cap (N'\cap N_\infty)^\omega=R'\cap
S^\omega$. Thus, if $\sigma \in \Cal G_0$ then $\psi \in \text{\rm
Ctr}(N_\infty)$.

Conversely, if $\theta\in \Cal H$ then define $\Phi(\theta)$ to be
the restriction of $\theta$ to $R=P'\cap N_\infty$. Note that by
Lemma 2.5 (a), we have $\Phi(\theta)\in \Cal G$. Also, if $\theta\in
\Cal H_0$ then $\Phi(\theta)\in \Cal G_0$. Since we clearly have
$\Phi \circ \Psi = id$ and $\Psi \circ \Phi=id$, the statement
follows.

\hfill $\square$

Note that, for a special class of inclusions of factors $P\subset
N$, the equality $\chi(N,P)=\chi(N_\infty)$ was already established
in [R95].

\proclaim{4.3. Theorem}  Let $P\subset N$ be a subfactor of finite
index. Assume $N$ has spectral gap in $N_\infty$ and $N'\cap
N_\infty$ is a factor $($e.g. $N$ is non-Gamma and $P\subset N$ has
finite depth$)$.

\vskip .05in \noindent $(i)$ $N_\infty$ is s-McDuff if and only if $P
\subset N$ is a ``matricial'' inclusion, i.e.
it is of the form $P \subset M_n(P)$, for some $n$, or equivalently
$N=P\vee (P'\cap N)$.

\vskip .05in \noindent $(ii)$ If $P'\cap N=\Bbb C$ and
$\Gamma_{P,N}$ is strongly amenable and its canonical weight vectors
$\vec{s}=(s_k)_k$, resp. $\vec{t}=(t_l)_l$, have distinct entries
$($i.e. $s_k\neq s_{k'}$ for $k\neq k'$ and $t_l\neq t_{l'}$ if
$l\neq l')$, then $\chi(N_\infty)=1$.
\endproclaim
\vskip .05in \noindent {\it Proof.} To prove $(i)$, assume
$Q\overline{\otimes} R = N_\infty$. Let $R_n\nearrow R$ be an
increasing sequence of finite dimensional subfactors exhausting $R$.
By Lemma 3.2 there exists $n$, such that $N\subset_{1/4}
Q\overline{\otimes} R_n$. Since $N'\cap N_\infty$ is a factor, by
(Proposition 12 in [OP03]), it follows that there exists a unitary
element $u\in N_\infty$ such that $uNu^*\subset Q^t$, for some
$t>0$. Thus, by replacing $Q$ by $u^*Q^tu$ and $R$ by $u^*R^{1/t}u$,
we may assume the decomposition $N_\infty=Q\overline{\otimes} R$ is
so that $N\subset Q$. Hence, $R=Q'\cap N_\infty \subset N'\cap
N_\infty$.

Now note that since $R$ splits off $N_\infty$, it also splits off
$N'\cap N_\infty$, i.e., if we denote $B=R'\cap (N'\cap N_\infty)$
then $N'\cap N_\infty = R \overline{\otimes} B$. In particular, $B$
is a factor.

On the other hand, by applying Lemma 3.1 again, for any $\varepsilon
> 0$ there exists $n$ such that $Q\subset_{\varepsilon/2} N_n$.
Taking relative commutants, it follows that $N_n'\cap N_\infty
\subset_{\varepsilon} R$. Thus, with the notations in [P03], if we
take $\varepsilon < 1$ then we get $N_n'\cap N_\infty
\prec_{N_\infty} R$, and thus $R\overline{\otimes} B=N'\cap
N_\infty\prec_{N_\infty} R$ as well (because $N_n'\cap N_\infty
\subset N'\cap N_\infty$ has finite index). In particular, this
shows that $B$ must be a finite dimensional factor. Thus, by
replacing if necessary $R$ by $R\overline{\otimes} B$, we may
actually assume $N'\cap N_\infty = R$. So $N'\cap N_\infty$ splits
off $N_\infty$. Similarly $P'\cap N_\infty$ splits off $N_\infty$ as
well, implying that $N'\cap N_\infty$ splits off $P'\cap N_\infty$.
Thus, $N'\cap N_\infty \subset P'\cap N_\infty$ is matricial, in
particular extremal, which in turn implies $(P'\cap N_\infty)'\cap
N_\infty = P\subset N=(N'\cap N_\infty)'\cap N_\infty$ is matricial
as well.

To prove $(ii)$, note that by Theorem 4.1 we have $\chi(N_\infty)=
\chi(P'\cap N_\infty, N'\cap N_\infty)$. Note also that by (part
$4^\circ$ of Theorem 1.6 in [P92]), any automorphism $\theta$ of
$N'\cap N_\infty \subset P'\cap N_\infty$ is either implemented by a
unitary $u\in P'\cap N_\infty$ normalizing $N'\cap N_\infty$ or is
properly outer. If $u$ is not in $N'\cap N_\infty$, then the
canonical vector of $\Gamma_{N'\cap N_\infty, P'\cap N_\infty}$
would have an even vertex other than $*$ equal to 1. By strong
amenability, $P=(P'\cap N_\infty)'\cap N_\infty \subset (N'\cap
N_\infty)'\cap N_\infty = N$, which implies the same is true for
$\Gamma_{P,N}$ and its standard vector, contradicting the
hypothesis.

Thus, either $\theta=\text{\rm Ad}(u)$ for some $u\in N'\cap
N_\infty$ of $\theta$ is properly outer on the inclusions $N'\cap
N_\infty\subset P'\cap N_\infty$. But by (part $6^\circ$ of Theorem
1.6 in [P92]), any properly outer automorphism of $N'\cap N_\infty
\subset P'\cap N_\infty$ is centrally free, i.e. it acts freely on
$N'\cap (P'\cap N_\infty)^\omega=N_\infty'\cap N_\infty^\omega$, a
contradiction. Thus, $\chi(P'\cap N_\infty, N'\cap N_\infty)$ must
be trivial, and so $\chi(N_\infty)=1$ as well. \hfill $\square$

\proclaim{4.4. Corollary}  For each $2\leq n \leq \infty$, the
factor $L(\Bbb F_n)\overline{\otimes} R$ has a sequence of
irreducible subfactors $M_m$ with distinct indices, satisfying
$\chi(M_m)=1$ and which are McDuff but not s-McDuff.
\endproclaim
\vskip .05in \noindent {\it Proof.} By [R92], for each $m\geq 5$
there exists a subfactor $P\subset N=L(\Bbb F_n)$ of index $4\cos^2
\pi/m$ and standard graph equal to $A_{m-1}$. Moreover, subfactors
with such graph have weight vectors with distinct entries,
$\vec{s_m}=(s_{m,k})_k$ and the numbers $r_m=\Sigma_k s_{m,k}^2$ are
distinct, for $m\geq 5$ (see e.g. [GHJ89] or [EK98]). Also, if we
fix $m$ and denote by $N_\infty$ the enveloping factor of the
inclusion $P\subset N$ as before and let $R=N'\cap N_\infty$, then
the inclusion of factors $L(\Bbb F_n)\vee R \subset N_\infty$ has
index equal to $r_m$ (see e.g. [EK98]). By [J82], there exist a
subfactor $M_m\subset L(\Bbb F_n)\vee R$, such that $L(\Bbb F_n)\vee
R \subset N_\infty$ is the basic construction of $M_m \subset L(\Bbb
F_n)\vee R$ and by [PP83], $M_m$ is an amplification of $N_\infty$.
Thus, since by Theorem 4.3 $N_\infty$ is not s-McDuff and
$\chi(N_\infty)=1$, it follows that $M_m$ satisfies the same
properties as well. \hfill $\square$

\proclaim{4.5. Theorem}  Let $P\subset N$ be an inclusion of
non-Gamma $\text{\rm II}_1$ factors with Temperley-Lieb standard
invariant, i.e., either $[N:P]< 4$ and $\Gamma_{P,N}=A_n$ for some
$n\geq 2$, or $[N:P]\geq 4$ and $\Gamma_{P,N}=A_\infty$. Assume $N$
has spectral gap in $N_\infty$ $($note that this is automatic in
case $[N:P] < 4)$. If $[N:P]\leq 4$ then $\chi(N_\infty)=1$, while
if $[N:P]>4$, then $\chi(N_\infty)=\Bbb T$.
\endproclaim
\vskip .05in \noindent {\it Proof.} The case $[N:P]\leq 4$ is a
consequence of 4.3.$(ii)$. To prove the case $[N:P]>4$, note that by
[P90], the inclusion of factors $(P'\cap N_\infty)'\cap
N_\infty=\tilde{P}\subset \tilde{N}=(N'\cap N_\infty)'\cap N_\infty$
is locally trivial. More precisely, there is an isomorphism $\sigma:
p\tilde{N}p \simeq (1-p)\tilde{N}(1-p)$, where $p\in \tilde{P}'\cap
\tilde{N}$ is a projection satisfying $\tau(p)\tau(1-p)=
[N:P]^{-1}$, such that $\tilde{P}=\{x + \sigma(x) \mid x\in
p\tilde{N}p\}$. Also, the enveloping factor of $\tilde{P}\subset
\tilde{N}$ is equal to $N_\infty$ and $\tilde{N}$ has spectral gap
in $N_\infty$ (because $N$ does). But then Corollary 2.7.1$^\circ$
applies, showing that any element in $\overline{\text{\rm
Int}(N_\infty)} \cap \text{\rm Ctr}(N_\infty)/\text{\rm
Int}(N_\infty)$ is represented by an automorphism $\theta$ of
$N_\infty$ which acts trivially on $\tilde{N}$ and for which there
exists a non-zero partial isometry $w\in \tilde{N}'\cap N_\infty$
and some $m$ such that $\theta(x)w=wx$, $\forall x\in
\tilde{N}_m'\cap N_\infty$. Note that $w^*w$ lies in
$(\tilde{N}_m'\cap N_\infty)'\cap (\tilde{N}'\cap
N_\infty)=\tilde{N}'\cap \tilde{N}_m$. By multiplying from the right
on both sides with a minimal projection $p\leq w^*w$ in
$\tilde{N}'\cap \tilde{N}_m$, we may in fact assume $w^*w=p$. Since
the minimal projections in distinct direct summands of
$\tilde{N}'\cap \tilde{N}_m$ have distinct traces, it follows that
both $ww^*$  and $\theta(p)$ are minimal projection in the same
direct summand of $\theta(\tilde{N}'\cap \tilde{N}_m)$. Thus, by
multiplying $w$ from the left with an appropriate partial isometry
in $\theta(\tilde{N}'\cap \tilde{N}_m)$, we may assume
$ww^*=\theta(p)$. So, now we have $\theta(x)w=wx$, $\forall x\in
(\tilde{N}_m'\cap N_\infty)p=p(\tilde{N}'\cap N_\infty)p$. But this
implies there exists a unitary element $v\in \tilde{N}'\cap
N_\infty$ such that $\theta=\text{\rm Ad}v$ on $\tilde{N}'\cap
N_\infty$.

Altogether, this shows that any element in $\chi(N_\infty)$ can be
represented by an automorphism acting trivially on $Q=\tilde{N}\vee
\tilde{N}'\cap N_\infty$, i.e. $\chi(N_\infty)$ naturally embeds
into the group $G$ of automorphisms of $N_\infty$ acting trivially
on $Q$.

But $N_\infty$ is isomorphic to the crossed product of $Q$ by an
aperiodic automorphism, whose amplification to $Q\overline{\otimes}
\Cal B(\Cal H) \overline{\otimes} \Cal B(\Cal H)$ acts by $\sigma
\otimes \sigma_0$, where $\sigma$ is an automorphism of
$\tilde{N}\overline{\otimes} \Cal B(\Cal H)$ scaling the trace by
$t/(1-t)$ and $\sigma_0$ is an automorphism of $\tilde{N}'\cap
N_\infty \overline{\otimes} \Cal B(\Cal H)\simeq R\overline{\otimes}
\Cal B(\Cal H)$ scaling the trace by $(1-t)/t$ (both given by the
inclusion $\tilde{N}'\cap N_\infty \subset \tilde{P}'\cap
N_\infty$). Thus $G=\Bbb T$, with each $\lambda \in \Bbb T$
corresponding to the automorphism $\theta_\lambda$ of $N_\infty = Q
\rtimes \Bbb Z$ which leaves $Q$ pointwise fixed and satisfies
$\theta_\lambda(u)=\lambda u$, where $u\in N_\infty$ is the
canonical unitary implementing the action of $\Bbb Z$ on $Q$.

To see that each $\theta_\lambda$ lies in $\overline{\text{\rm
Int}(N_\infty)} \cap \text{\rm Ctr}(N_\infty)$, note that one can
take a sequence of unitary elements $u_n$ in the hyperfinite II$_1$
factor $R=\tilde{N}'\cap N_\infty$ such that $\sigma_0(u_n)=\lambda
u_n$ and $(u_n)_n \in R'\cap R^\omega$. But then $\theta=\lim_n
\text{\rm Ad}u_n$ acts trivially on $Q$ while $\theta(u)=\lambda u$.
Moreover, $\theta$ clearly belongs to $\overline{\text{\rm
Int}(N_\infty)} \cap \text{\rm Ctr}(N_\infty)$, by the way it was
defined. \hfill $\square$

\head  References \endhead

\item{[C75]} A. Connes: {\it Sur la classification des facteurs
de type} II, C. R. Acad. Sci. Paris {\bf 281} (1975), 13-15.

\item{[C76]} A. Connes: {\it Classification of injective factors}, Ann.
Math. {\bf 104} (1976),  73-115.

\item{[EK98]} D. E. Evans, Y. Kawahigashi: ``Quantum symmetries on
operator algebras'', Oxford University Press, 1998.

\item{[GHJ90]} F. Goodman, P. de la Harpe, V.F.R. Jones: ``Coxeter graphs and
towers of algebras'', MSRI Publications {\bf 14}, Springer, 1989.

\item{[J82]} V.F.R. Jones: {\it Index for subfactors},
Invent. Math {\bf 72} (1983), 1-25.

\item{[K93]} Y. Kawahigashi: {\it Centrally trivial automorphisms and
an analogue of Connes's} $\chi(M)$ {\it for subfactors}, Duke Math.
J. {\bf 71} (1993), 93-118.

\item{[M70]} D. McDuff: {\it Central sequences and the hyperfinite factor}, Proc.
London Math. Soc. {\bf 21} (1970), 443–461.

\item{[Oc87]} A. Ocneanu: {\it Quantized groups, string algebras and Galois
theory for von Neumann algebras}. In ``Operator Algebras and
Applications'', London Math. Soc. Lect. Notes Series, Vol. {\bf
136}, 1988, pp. 119-172.

\item{[OP03]} N. Ozawa, S. Popa: {\it Some prime factorization results for type}
II$_1$ {\it factors}, Invent. Math., {\bf 156} (2004), 223-234.

\item{[OP07]} N. Ozawa, S. Popa: {\it On a class of $\text{\rm II}_1$
factors with at most one Cartan subalgebra}, math.OA/0706.3623, to
appear in Annals of Mathematics

\item{[PP83]} M. Pimsner, S. Popa: {\it Entropy and index for subfactors}, Ann.
Sc. Ec. Norm. Sup. {\bf 19} (1986), 57-106.

\item{[P90]} S. Popa: {\it Markov traces on universal
Jones algebras and subfactors of finite index}, Invent. Math. {\bf
111} (1993), 375-405.

\item{[P92]} S. Popa: {\it Classification of actions of
discrete amenable groups on amenable subfactors of type} II, IHES
preprint 1992, to appear in Intern. J. Math, 2009
(see http://www.math.ucla.edu/~popa/preprints.html)

\item{[P94]} S. Popa, {\it An axiomatization of the lattice of
higher relative commutants of a subfactor}, Invent. Math., {\bf 120}
(1995), 427-445.

\item{[98]} S. Popa: {\it Universal construction of subfactors}, J.
reine angew. Math., {\bf 543} (2002), 39-81.

\item{[P01]} S. Popa: {\it On a class of type} II$_1$ {\it factors with Betti
numbers invariants}, Ann. of Math {\bf 163} (2006), 809-899.

\item{[P03]} S. Popa: {\it Strong Rigidity of} II$_1$ {\it Factors
Arising from Malleable Actions of $w$-Rigid Groups} I, Invent. Math.
{\bf 165} (2006), 369-408 (math.OA/0305306).

\item{[P06]} S. Popa: {\it On the superrigidity
of malleable actions with spectral gap}, J. Amer. Math. Soc. {\bf
21} (2008), 981-1000 (math.GR/0608429).

\item{[PS00]} S. Popa, D. Shlyakhtenko: {\it Universal properties of}
$L(\Bbb F_\infty)$ {\it in subfactor theory}, Acta Mathematica, {\bf
191} (2003), 225-257.

\item{[R92]} F. Radulescu: {\it Random matrices, amalgamated free products and
  subfactors of the von Neumann algebra of a free group, of noninteger
  index}, Invent. Math. {\bf 115} (1994), 347-389.

\item{[R95]} F. Radulescu: {\it A new invariant for subfactors of the
von Neumann algebra of a free group}, Proceedings of The Fields
Institute Workshop in Free Probability, D. Voiculescu ed, 1995, pp
213-240.

\enddocument